\documentclass[11pt]{article}
\usepackage[brazil,english]{babel}
\usepackage[latin1]{inputenc}
\usepackage{fullpage}
\usepackage[dvips]{graphicx}
\usepackage{amssymb}
\usepackage{amsmath}
\usepackage{color}
\usepackage{footnote}
\usepackage{hyperref}

\newtheorem{theorem}{Theorem}[section] 
\newtheorem{corollary}[theorem]{Corollary}  
\newtheorem{proposition}[theorem]{Proposition}  
\newtheorem{lemma}[theorem]{Lemma}  
\newtheorem{definition}[theorem]{Definition} 
\newtheorem{example}[theorem]{Example}

\def\beginproof{\noindent{\em Proof}.~}
\def\endproof{{\ \hfill\hbox{\fbox{}}\parfillskip 0pt}\par}

\renewcommand{\thefootnote}{\arabic{ftnote}}

\newcommand{\R}{\mathbb{R}}

\definecolor{mael}{RGB}{250, 0, 250} 
\definecolor{ad}{RGB}{50, 0, 250}
\definecolor{cmt}{RGB}{0, 150, 20}

\definecolor{rmv}{RGB}{250, 0, 0}
\title{On strong second-order optimality conditions under \\ relaxed constant rank 
constraint qualification}

\author{Ademir A. Ribeiro\footnotemark[4] \and Mael Sachine\footnotemark[4]}

\begin{document}

\maketitle
\renewcommand{\thefootnote}{\fnsymbol{footnote}}
\footnotetext[4]{Department of Mathematics, Federal University of Paran\'a, 
Brazil (ademir.ribeiro@ufpr.br, mael@ufpr.br).}

\renewcommand{\thefootnote}{\arabic{ftnote}}

\begin{abstract}
We discuss the (first- and second-order) optimality conditions for nonlinear programming 
under the relaxed constant rank constraint qualification. This condition generalizes 
the so-called linear independence constraint qualification. Although the optimality 
conditions are well established in the literature, the proofs presented here are based 
solely on the well-known inverse function theorem. This is the only prerequisite from 
real analysis used to establish two auxiliary results needed to prove the optimality 
conditions, thereby making this paper totally self-contained.
\end{abstract}


{\bf Subclass. }
90C30, 90C33, 90C46 

\thispagestyle{plain}

\section{Introduction}
In this short note we are concerned with the nonlinear programming problem 
\begin{equation}
\label{problem}
\begin{array}{cl}
\displaystyle\mathop{\rm minimize }  & f(x)  \\
{\rm subject\ to } & g(x)\leq 0, \\
& h(x)=0, 
\end{array}
\end{equation}
where $f:\R^n\to\R$, $g:\R^n\to\R^m$ and $h:\R^n\to\R^p$ are twice continuously 
differentiable functions. The feasible set of the problem (\ref{problem}) is denoted by 
$$
\Omega=\{x\in\R^n\mid g(x)\leq 0, h(x)=0\}
$$
and the active index set by $I_g(x)=\{i\mid g_i(x)=0\}$. Also, the Lagrangian function 
associated with \eqref{problem} is $\ell:\R^n\times\R^m\times\R^p\to\R$ given by
\begin{equation}
\label{Lagrangian}
\ell(x,\mu,\lambda)=
f(x)+\mu^Tg(x)+\lambda^Th(x).
\end{equation}

The main role of the several constraint qualifications (CQs) 
\cite{Abadie67,AndreaniHaeserSchuverdtSilva12b,AndreaniMartinezRamosSilva16,AndreaniMartinezSchuverdt,Guignard69,Janin84,MangasarianFromovitz,MinchenkoStakhovski11a,QiWei} 
is to guarantee the existence of Lagrange multipliers for any continuously differentiable 
function $f$ at the points of its local minimum on the set $\Omega$ or in other words the 
satisfaction of the Karush-Kuhn-Tucker conditions 
$$
\nabla_{x}\ell(x,\mu,\lambda)=0,\quad\mu^Tg(x)=0,\quad\mu\geq 0.
$$
Moreover, most textbooks on nonlinear programming present the classical second-order result: 
under linear independence constraint qualification (LICQ) the Lagrangian Hessian is positive 
semidefinite on the so-called {\em critical cone}. 

Nevertheless, not all constraint qualifications have this obviously desirable property. 
Indeed, many of the most important CQs do not ensure such a positive semidefiniteness 
(see, for instance, Example \ref{ex:rcrcq}). 

Here we focus our attention in two second-order CQs that are natural generalizations of LICQ. 

\begin{definition}[\cite{Janin84}]
\label{def:crcq}
We say that the {\em constant rank constraint qualification} (CRCQ) holds at a feasible point 
$\bar{x}\in \Omega$ if for any subsets $\mathcal I\subset I_g(\bar{x})$ and 
$\mathcal J\subset \{1,\ldots,p\}$, the family of gradients
$$
\nabla g_i(x), \ i\in \mathcal{I}, \quad \nabla h_j(x), \ j\in \mathcal{J}, 
$$
has the same rank for every $x$ in a neighborhood of $\bar{x}$.
\end{definition}

\begin{definition}[\cite{MinchenkoStakhovski11a}]
\label{def:rcrcq}
We say that the {\em relaxed constant rank constraint qualification} (RCRCQ) holds at a 
feasible point $\bar{x}\in \Omega$ if for any subset $\mathcal I\subset I_g(\bar{x})$, the 
family of gradients
$$
\nabla g_i(x), \ i\in \mathcal I, \quad \nabla h_j(x), \ j\in\{1,\ldots,p\}, 
$$
has the same rank for every $x$ in a neighborhood of $\bar{x}$.
\end{definition}

Note that the only difference between the above definitions relies on the set of equality 
gradients, so that for inequality constrained problems they are equivalent. 
We mention that there are equivalent ways to characterize RCRCQ and CRCQ in terms 
of linear dependence of vectors \cite{AndreaniHaeserSchuverdtSilva12b,QiWei}. It is also 
obvious that these two conditions generalize LICQ. On the other hand, (R)CRCQ and 
Mangasarian-Fromovitz constraint qualification (MFCQ) are independent of 
each other \cite{AndreaniHaeserSchuverdtSilva12b,Janin84}.

In \cite{AndreaniEchagueSchuverdt10}, the authors prove strong second-order necessary 
conditions (SSONC) under CRCQ. Assuming CRCQ, it is proved in \cite{Janin84} the 
first-order Kuhn-Tucker regularity condition. In \cite[Theorem 1]{MinchenkoStakhovski11a} 
and \cite[Lemma 6]{MinchenkoStakhovski11b} it is proved that RCRCQ implies Abadie 
Constraint Qualification (ACQ) by making use of Lyusternik's theorem and rank theorem, 
respectively. In \cite[Corollary 2.3]{MinchenkoLeschov16} and 
\cite[Theorem 6]{MinchenkoStakhovski11b} it is established SSONC assuming RCRCQ. 
Other assumptions ensuring second-order necessary conditions can be found in
\cite{AndreaniBehlingHaeserSilva17,AndreaniHaeserRamosSilva17}, as we will detail later.

In this paper we revisit first- and second-order optimality conditions associated with 
the problem (\ref{problem}). Concerning the prerequisites for the main results of the paper, 
we generalize some results of \cite{AndreaniEchagueSchuverdt10,Janin84} assuming 
RCRCQ instead of CRCQ. We then establish that RCRCQ is indeed a constraint qualification and 
that SSONC holds by means of a different, simpler and self-contained analysis (when 
comparing with \cite{MinchenkoLeschov16,MinchenkoStakhovski11a,MinchenkoStakhovski11b}), 
using only the inverse function theorem as the basic tool to obtain auxiliary results.

The paper is organized as follows: in Sect. \ref{sec:prelim} we provide the basic mathematical 
tools used in this paper. Sect. \ref{sec:rcrcq} is devoted to discuss the first- and 
second-order optimality conditions under RCRCQ. Conclusions are presented in 
Sect. \ref{sec:concl}.

{\noindent\bf Notation.} Throughout this paper, the full Jacobian (or derivative) of a 
function $\xi:\R^n\to\R^s$ at $x\in\R^n$ is denoted by 
$\xi'(x)=\nabla\xi(x)^T\in\R^{s\times n}$, while for an 
index set $J\subset\{1,\ldots,n\}$ such that $|J|=r$ the partial Jacobian is 
$\partial_J\xi(x)=\nabla_{x_J}\xi(x)^T\in\R^{s\times r}$.

\section{Preliminaries}
\label{sec:prelim}
In this section we prove two technical lemmas based on mathematical analysis which will be 
used to establish the main results of this paper. 

The first auxiliary lemma says that a function of rank $r$ can be expressed (after a 
change of variables) in terms of $r$ variables in a very special way. This result is 
established based solely on the well-known inverse function theorem.

\begin{lemma}
\label{lm:rank1}
Let $\sigma:\R^n\to\R^q$ be a twice continuously differentiable function and suppose that 
$\sigma'(x)$ has constant rank, say $r$, in a neighborhood of the point $\bar{x}\in\R^n$. 
Then there exist an open set $Z\subset\R^n$, a neighborhood $V$ of $\bar{x}$, an index set 
$J\subset\{1,\ldots,n\}$ and twice continuously differentiable functions $\alpha:Z\to V$ and 
$\lambda:\R^r\to\R^{q-r}$ such that $|J|=r$ and $\sigma(\alpha(z))=(z_J,\lambda(z_J))$ for 
all $z\in Z$.
\end{lemma}
\beginproof
Suppose, after rearranging the components of $\sigma$ if necessary, that 
$${\rm rank}\{\nabla\sigma_1(\bar{x}),\ldots,\nabla\sigma_r(\bar{x})\}=r$$ and define 
$\xi:\R^n\to\R^r$ by $\xi(x)=(\sigma_1(x),\ldots,\sigma_r(x))$. Then there exists an index set 
$J\subset\{1,\ldots,n\}$ such that $|J|=r$ and $\partial_J\xi(\bar{x})\in\R^{r\times r}$ is 
nonsingular. Define $K=\{1,\ldots,n\}\setminus J$ and $\varphi:\R^n\to\R^n$ by 
$\varphi(x)=(\xi(x),x_K)$. Thus, 
$$
\varphi'(\bar{x})=\left(\hspace{-2pt}\begin{array}{cc} 
\partial_J\xi(\bar{x}) & \partial_K\xi(\bar{x}) \\ 0 & I \end{array}\hspace{-2pt}\right)
$$
is nonsingular and hence we can apply the inverse function theorem to conclude that there exist 
a neighborhood $V$ of $\bar{x}$ and an open set $Z\subset\R^n$ such that $\varphi_{\mid V}:V\to Z$ 
is invertible with inverse $\alpha:Z\to V$ twice continuously differentiable. 
Writing $\alpha(z)=(\alpha_J(z),\alpha_K(z))$, for all $z\in Z$ we have 
$$
z=\varphi(\alpha(z))=(\xi(\alpha(z)),\alpha_K(z)),
$$
which means that 
\begin{equation} 
\label{eq1_lm:rank}
\xi(\alpha(z))=z_J\quad\mbox{and}\quad\alpha_K(z)=z_K.
\end{equation}
Therefore, letting $\lambda(z)=(\sigma_{r+1}(\alpha(z)),\ldots,\sigma_q(\alpha(z)))$, we have 
$$
\sigma(\alpha(z))=(\xi(\alpha(z)),\lambda(z))=(z_J,\lambda(z)).
$$
It remains to show that $\lambda(z)$ does not depend on $z_K$. For this purpose note that 
$$
\sigma'(\alpha(z))\alpha'(z)=(\sigma\circ\alpha)'(z)=\left(\hspace{-2pt}\begin{array}{cc}
I & 0 \\ 
\partial_J\lambda(z) & \partial_K\lambda(z) \end{array}\hspace{-2pt}\right).
$$
Moreover, since we may assume that $\sigma'$ has constant rank $r$ in $V$, we have 
${\rm rank}(\sigma'(\alpha(z))\alpha'(z))=r$, which implies that $\partial_K\lambda(z)=0$.
Moreover, the set $Z$ can be assumed to be convex (otherwise, take an open ball 
$B(\varphi(\bar{x}),\delta)\subset Z$ and redefine the sets: 
$Z\leftarrow B(\varphi(\bar{x}),\delta)$ and $V\leftarrow\alpha(Z)$).
Therefore, $\lambda(z)$ does not depend on $z_K$ and can be written as $\lambda(z_J)$. 
\endproof

\begin{figure}[htbp]
\centering
\includegraphics[scale=.42]{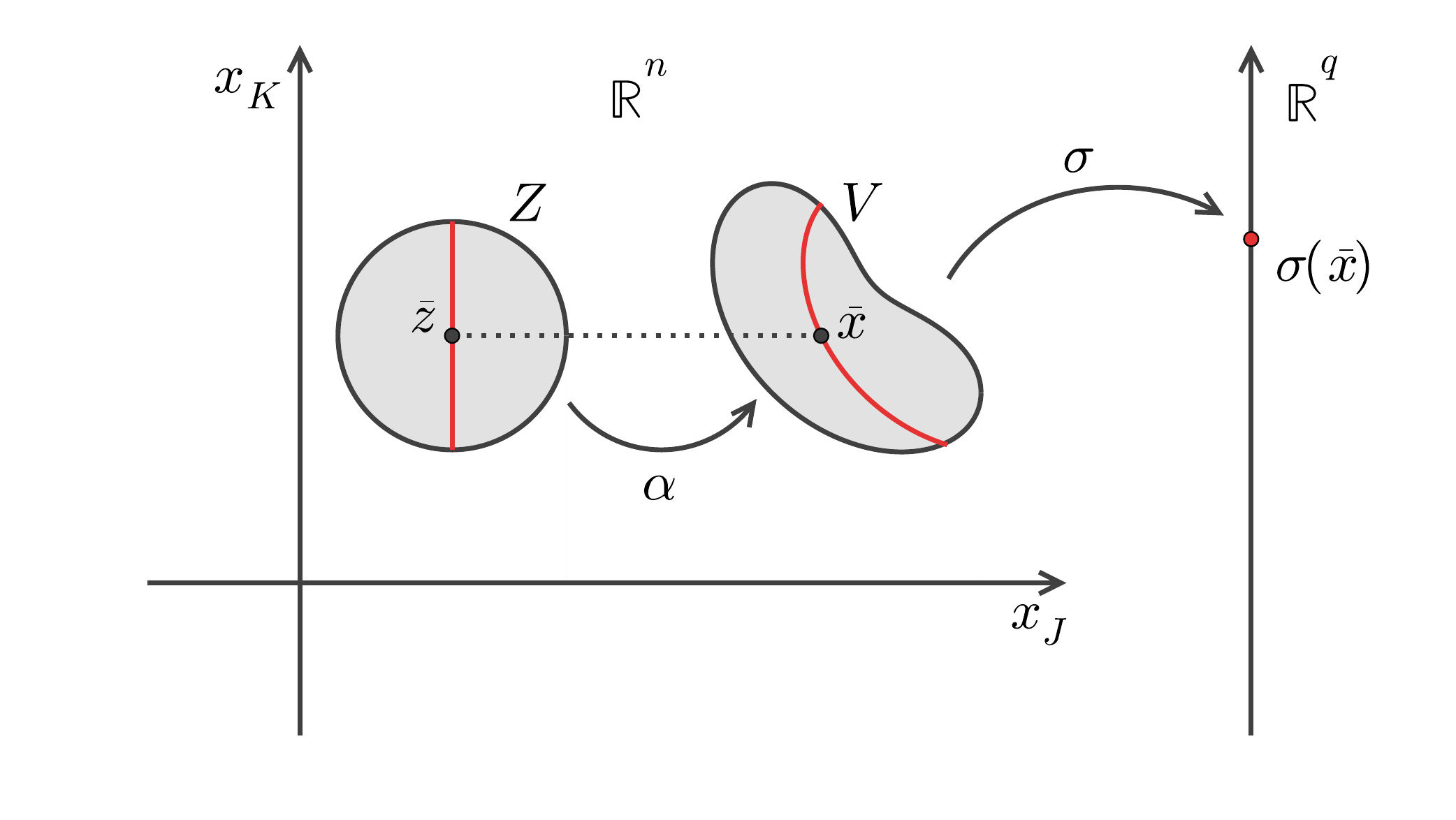}
\caption{Graphical illustration of Lemma \ref{lm:rank1}. The function $\alpha$ preserves 
the $K$ coordinates: $\alpha(z)=(\alpha_J(z),z_K)$. On the other hand, since 
$\sigma(\alpha(z))=(z_J,\lambda(z_J))$, points on the same ``vertical line'' (constant 
$J$ coordinates) are mapped by $\sigma\circ\alpha$ into a single point in $\R^q$.}
\label{fig1}
\end{figure}

Figure \ref{fig1} illustrates Lemma \ref{lm:rank1}. We point out that $z$ and $\alpha(z)$ 
share the same $K$ coordinates, as established in the relation (\ref{eq1_lm:rank}) of the 
proof. Moreover, since $\sigma(\alpha(z))=(z_J,\lambda(z_J))$, we have that all the points 
on the same ``vertical line'' (constant $J$ coordinates) are mapped by $\sigma\circ\alpha$ 
into a single point in $\R^q$.

In the next result we rewrite \cite[Proposition 3.1]{AndreaniEchagueSchuverdt10} 
and \cite[Proposition 2.2]{Janin84} in a suitable and more precise way. Despite the authors 
in these references comment, without adding any further argument, that this result is a special 
case of the constant rank theorem of \cite{Malliavin72}, we consider that such a relation 
is not evident and the analysis is somewhat involved. We then give here a completely 
self-contained proof, using only the previous lemma.

\begin{lemma}
\label{lm:rank2}
Let $\sigma:\R^n\to\R^q$ be a twice continuously differentiable function and suppose that 
$\sigma'(x)$ has constant rank in a neighborhood of the point $\bar{x}\in\R^n$. Denote 
$D=\mathcal{N}(\sigma'(\bar{x}))$ the nullspace of $\sigma'(x)$. Then there exist 
neighborhoods $U$, $V$ of $\bar{x}$ and a twice continuously differentiable function 
$\phi:U\to V$ such that 
\begin{subequations}
\begin{align}
\phi(\bar{x})=\bar{x}, \label{phi1} \\ 
\phi'(\bar{x})=I\in\R^{n\times n} \, (\mbox{the identity matrix}), \label{phi2} \\ 
\sigma(\phi(x+d))=\sigma(\phi(x)) 
\mbox{ for all $x\in U$ and $d\in D$ such that $x+d\in U$.} \label{phi3}
\end{align}
\end{subequations}
In particular, $\sigma(\phi(\bar{x}+d))=\sigma(\bar{x})$ for all $d\in D$ such that 
$\bar{x}+d\in U$.
\end{lemma}
\beginproof
Let ${\rm rank}(\sigma'(\bar{x}))=r$ and consider the function $\alpha:Z\to V$ defined in 
Lemma \ref{lm:rank1}, which we know that has the form $\alpha(z)=(\alpha_J(z),z_K)$. Define 
$\bar{z}=\alpha^{-1}(\bar{x})=(\xi(\bar{x}),\bar{x}_K)$, 
$A=(\alpha'(\bar{z}))^{-1}\in\R^{n\times n}$, $\gamma:\R^n\to\R^n$ by 
$$
\gamma(x)=A(x-\bar{x})+\bar{z},
$$
$U=\gamma^{-1}(Z)$ and $\phi=\alpha\circ\gamma:U\to V$. Thus, we note immediately that 
$\phi(\bar{x})=\bar{x}$, giving (\ref{phi1}). Moreover, (\ref{phi2}) follows from 
$$
\phi'(\bar{x})=\alpha'(\gamma(\bar{x}))\gamma'(\bar{x})=\alpha'(\bar{z})\gamma'(\bar{x})=I.
$$
Finally, take $x\in U$ and $d\in D$ such that $x+d\in U$. Then, 
$$
\sigma'(\alpha(\bar{z}))\alpha'(\bar{z})Ad=
\sigma'(\bar{x})\alpha'(\bar{z})Ad=\sigma'(\bar{x})d=0.
$$
On the other hand, since $\sigma(\alpha(z))=(z_J,\lambda(z_J))$, we have 
$$
\sigma'(\alpha(\bar{z}))\alpha'(\bar{z})=\left(\begin{array}{cc} 
I & 0 \\ \lambda'(\bar{z}_J) & 0 \end{array}\right).
$$
Therefore, 
$$
\left(\begin{array}{cc} I & 0 \\ \lambda'(\bar{z}_J) & 0 \end{array}\right)
\left(\begin{array}{c} {[Ad]}_J \\ {[Ad]}_K \end{array}\right)=
\left(\begin{array}{c} 0 \\ 0 \end{array}\right),
$$
which implies $[Ad]_J=0$ and hence, 
$$
\begin{array}{rcl}
\sigma(\phi(x+d)) & = & \sigma(\alpha(\gamma(x+d))) \vspace{.1cm} \\ 
& = & \sigma(\alpha(A(x-\bar{x}+d)+\bar{z})) \vspace{.1cm} \\ 
& = & ([A(x-\bar{x}+d)+\bar{z}]_J,\lambda([A(x-\bar{x}+d)+\bar{z}]_J)) \vspace{.1cm} \\ 
& = & ([A(x-\bar{x})+\bar{z}]_J,\lambda([A(x-\bar{x})+\bar{z}]_J)) \vspace{.1cm} \\ 
& = & ([\gamma(x)]_J,\lambda([\gamma(x)]_J)) \vspace{.1cm} \\ 
& = & \sigma(\alpha(\gamma(x))) \vspace{.1cm} \\ 
& = & \sigma(\phi(x)),
\end{array}
$$
proving (\ref{phi3}).
\endproof

\begin{figure}[htbp]
\centering
\includegraphics[scale=.42]{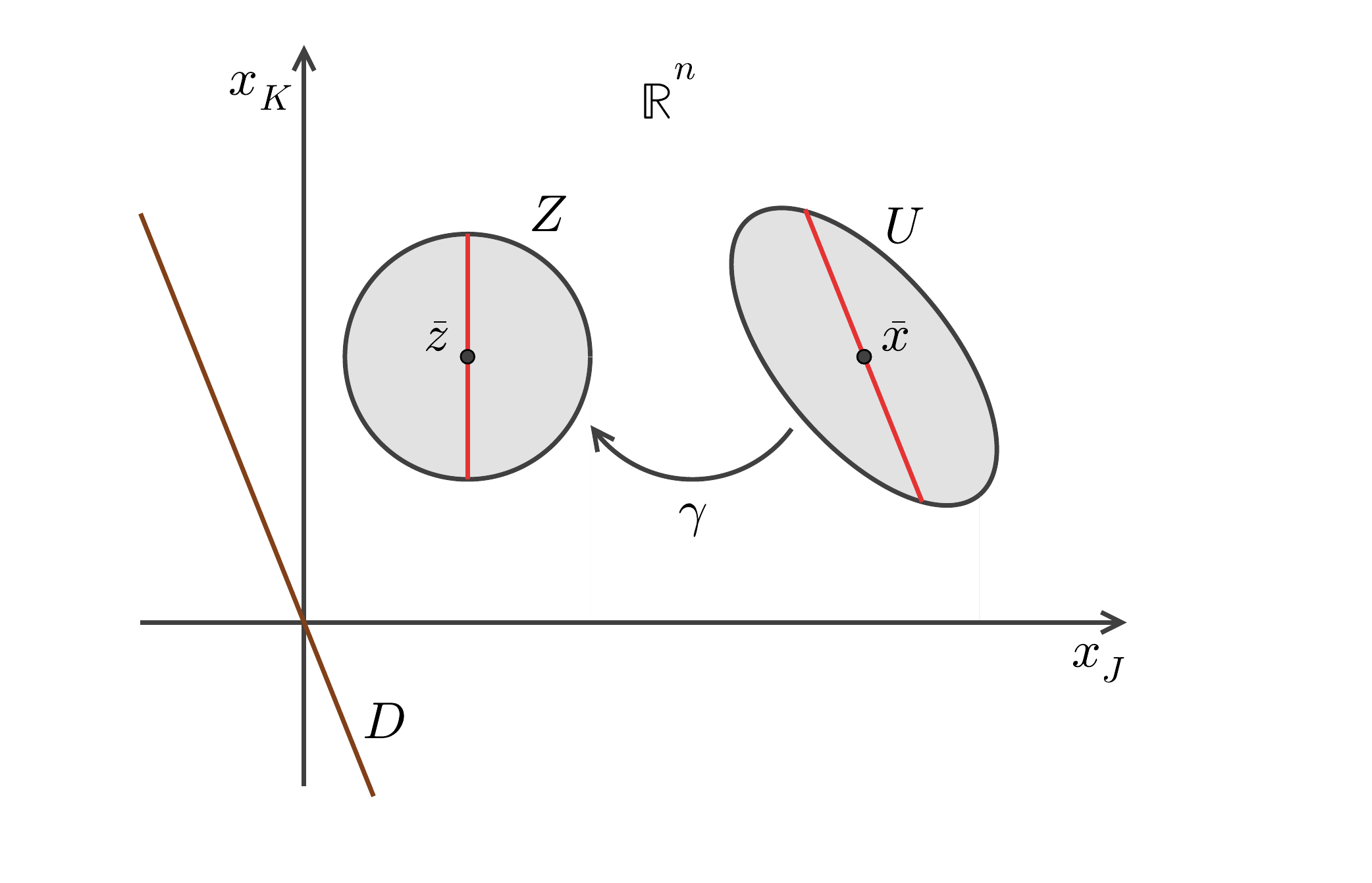}
\caption{Illustration of Lemma \ref{lm:rank2}. The linear function $\gamma$ maps 
the points $x+d\in U$ into the ``vertical line'' in $Z$, making constant the function $d\mapsto\sigma(\phi(x+d))$.}
\label{fig2}
\end{figure}

Lemma \ref{lm:rank2} is illustrated in Figure \ref{fig2}. The main role of the linear 
function $\gamma$ consists in mapping the points $x+d\in U$ into the ``vertical line'' 
in $Z$, making constant the function $d\mapsto\sigma(\phi(x+d))$. 

\section{Second-order optimality conditions}
\label{sec:rcrcq}
In order to establish the results of this section, let us recall the three cones 
associated with problem (\ref{problem}), 
the {\em tangent cone} to $\Omega$ at $\bar{x}\in\Omega$ 
$$
\mathcal{T}(\bar{x})=\left\{d\in\R^n\mid\exists\,(x^k)\subset\Omega\mbox{, } 
(t_k)\subset\R_+ \mbox{ : }t_k\to 0 \mbox{ and } 
\dfrac{x^k-\bar{x}}{t_k}\to d\right\},
$$
the {\em linearized cone} 
$$
\mathcal{L}(\bar{x})
=\{d\in\R^n\mid\nabla g_i(\bar{x})^Td\leq 0, \; i\in I_g(\bar{x}), \; 
\nabla h(\bar{x})^Td=0\}
$$
and the {\em strong critical cone} 
$$
\mathcal{C}^S(\bar{x})=
\{d\in\mathcal{L}(\bar{x}) \mid \nabla f(\bar{x})^Td \leq 0 \}.
$$

First of all, note that if $\bar{x}$ is a KKT point for problem \eqref{problem}, then 
we can rewrite $\mathcal{C}^S(\bar{x})$ without using $\nabla f$ but instead, using 
the multipliers of $\bar{x}$. To be precise, if 
$(\bar{x},\mu,\lambda)\in\R^n\times\R^m_+\times\R^p$ is a primal-dual KKT point and 
$I_g^+(\bar{x})=\{i\in I_g(\bar{x})\mid\mu_i>0\}$, then 
$$
\mathcal{C}^S(\bar{x})=
\{d\in\mathcal{L}(\bar{x})\mid\nabla g_i(\bar{x})^Td = 0, \ i\in I_g^+(\bar{x})\}.
$$
This follows immediately from the fact that 
$$
\nabla f(\bar{x})^Td=-\sum_{i\in I_g(\bar{x})} \mu_i \nabla g_i(\bar{x})^Td
$$
for all $d\in \mathcal{L}(\bar{x})$.

\begin{definition}
\label{def:ssonc}
Let $\bar{x}$ be a KKT point for problem \eqref{problem}. We say that $\bar{x}$ fulfills 
the {\em strong second-order optimality necessary condition} (SSONC) if for every 
multiplier vector $(\mu,\lambda)\in \R_+^m\times \R^p$ associated with $\bar{x}$ we have 
$$
d^T \nabla^2_{xx} \ell(\bar{x},\mu,\lambda) d\geq 0
$$
for all $d\in\mathcal{C}^S(\bar{x})$.
\end{definition}

In the example below, we see that MFCQ is not a second-order constraint qualification 
(and hence, none of CQs implied by MFCQ). 

\begin{example}
\label{ex:rcrcq}
Consider the problem in $\R^2$ given by 
$$
\begin{array}{cl}
\displaystyle\mathop{\rm minimize }  & f(x)=x_2  \\
{\rm subject\ to } & g_1(x)=x_1^2+(x_2-1)^2-1\leq 0, \\
& g_2(x)=1-x_1^2-(x_2+1)^2\leq 0.
\end{array}
$$
It is easy to see that $x^*=\left(\begin{array}{c} 0 \\ 0 \end{array}\right)$ is the global 
solution of the problem and that MFCQ holds at $x^*$. 
Indeed, $(x_2-1)^2\leq x_1^2+(x_2-1)^2\leq 1$ implies that $x_2\geq 0$. Moreover, we have 
$$
\nabla g_1(x^*)=\nabla g_2(x^*)=\left(\begin{array}{c} 0 \\ -2 \end{array}\right),
$$
which means that there exists $d\in\R^2$ such that $\nabla g_i(x^*)^Td<0$. 

On the other hand, 
$$
\nabla f(x^*)+\mu_1\nabla g_1(x^*)+\mu_2\nabla g_2(x^*)=
\left(\begin{array}{c} 0 \\ 1 \end{array}\right)+
\left(\begin{array}{c} 0 \\ -2\mu_1-2\mu_2 \end{array}\right),
$$
so that the Lagrange multipliers satisfy $2(\mu_1+\mu_2)=1$. The critical cone is 
$$
\mathcal{C}^S(x^*)=\{d\in\R^2\mid d_2=0 \}
$$
and the Lagrangian Hessian at this point is 
$$
\nabla^2f(x^*)+\mu_1\nabla^2g_1(x^*)+\mu_2\nabla^2g_2(x^*)=2(\mu_1-\mu_2)I.
$$
Thus, taking $\mu_1<\mu_2$ and $d\in\mathcal{C}^S(x^*)\setminus\{0\}$, we have 
$$
d^T\left(\nabla^2f(x^*)+\mu_1\nabla^2g_1(x^*)+\mu_2\nabla^2g_2(x^*)\right)d=
2(\mu_1-\mu_2)d_1^2<0.
$$
\end{example}

The result below shows, in particular, that (under RCRCQ) any linearized feasible 
direction $\bar{d}\in\mathcal{L}(\bar{x})$ is the velocity vector of a feasible arc at 
$\bar{x}\in\Omega$, which will allow us to derive ACQ directly from RCRCQ (without 
following the path of Figure 2 of \cite{AndreaniMartinezRamosSilva16}).

\begin{proposition}
\label{prop:arc}
Suppose that the relaxed constant rank constraint qualification holds at $\bar{x}\in\Omega$. 
Given $\bar{d}\in\mathcal{L}(\bar{x})$, define 
$J=\{j\in I_g(\bar{x})\mid\nabla g_j(\bar{x})^T\bar{d}=0\}$. Then there exists a twice 
continuously differentiable arc $\zeta:(-\delta,\delta)\to\R^n$ such that 
\begin{subequations}
\begin{align}
\zeta(0)=\bar{x}, \; \zeta'(0)=\bar{d}, \label{arc1} \\ 
g_j(\zeta(t))=0 \mbox{ for all $t\in(-\delta,\delta)$ and $j\in J$,} \label{arc2} \\ 
g_j(\zeta(t))\leq 0 \mbox{ for all $t\in(-\delta,\delta)$ and $j\notin I_g(\bar{x})$,} 
\label{arc3} \\ 
g_j(\zeta(t))\leq 0 \mbox{ for all $t\in[0,\delta)$ and $j\in I_g(\bar{x})\setminus J$,} 
\label{arc4} \\ 
h(\zeta(t))=0 \mbox{ for all $t\in(-\delta,\delta)$.} \label{arc5} 
\end{align}
\end{subequations}
In particular, $\zeta(t)\in\Omega$ for all $t\in[0,\delta)$ and hence, 
$\bar{d}\in\mathcal{T}(\bar{x})$.
\end{proposition}
\beginproof
Define $\sigma:\R^n\to\R^{|J|+p}$ 
defined by $\sigma(x)=(g_J(x),h(x))$. By the constraint qualification, $\sigma'(x)$ has constant 
rank in a neighborhood of $\bar{x}$. Therefore, applying Lemma \ref{lm:rank2} we conclude 
that there exist neighborhoods $U$, $V$ of $\bar{x}$ and a twice continuously differentiable 
function $\phi:U\to V$ such that $\phi(\bar{x})=\bar{x}$, $\phi'(\bar{x})=I\in\R^{n\times n}$ and 
\begin{equation}
\label{sigma_cte}
\sigma(\phi(\bar{x}+d))=\sigma(\bar{x})
\end{equation}
for all $d\in\mathcal{N}(\sigma'(\bar{x}))$ such that $\bar{x}+d\in U$. Consider $\delta>0$ 
such that $\bar{x}+t\bar{d}\in U$ for all $t\in(-\delta,\delta)$ and define 
$\zeta:(-\delta,\delta)\to\R^n$ by $\zeta(t)=\phi(\bar{x}+t\bar{d})$. Thus, 
$$
\zeta(0)=\phi(\bar{x})=\bar{x}\quad\mbox{and}\quad\zeta'(0)=\phi'(\bar{x})\bar{d}=\bar{d},
$$
giving (\ref{arc1}). Furthermore, since $\bar{d}\in\mathcal{N}(\sigma'(\bar{x}))$, we have 
$$
\sigma(\zeta(t))=\sigma(\phi(\bar{x}+t\bar{d}))=\sigma(\bar{x})=0,
$$
proving (\ref{arc2}) and (\ref{arc5}). The continuity of $g$ immediately gives (\ref{arc3}). 
Now, given $j\in I_g(\bar{x})\setminus J$, we have $\nabla g_j(\bar{x})^T\bar{d}<0$ and then 
$$
(g_j\circ\zeta)'(0)=g_j'(\zeta(0))\zeta'(0)=g_j'(\bar{x})\bar{d}=
\nabla g_j(\bar{x})^T\bar{d}<0,
$$
giving $g_j(\zeta(t))\leq g_j(\zeta(0))=0$ for all $t\geq 0$ sufficiently small. 
By reducing $\delta$ if necessary, we obtain (\ref{arc4}).
Finally, taking $(t_k)\subset(0,\delta)$ with $t_k\to 0$, we have 
$$
\dfrac{\zeta(t_k)-\bar{x}}{t_k}\to\bar{d},
$$
which means that $\bar{d}\in\mathcal{T}(\bar{x})$.
\endproof

We point out that Proposition \ref{prop:arc} generalizes 
\cite[Proposition 3.2]{AndreaniEchagueSchuverdt10} and \cite[Proposition 2.3]{Janin84} in 
three aspects. First, we assume the weaker hypothesis RCRCQ instead of CRCQ. Second, 
the arc $\zeta$ we constructed here is defined on the whole interval $(-\delta,\delta)$, 
while in the mentioned references the domain is an interval of the form $(0,\delta)$ and 
hence they have only the right-hand derivative. Finally we proved stronger smoothness of 
the arc $\zeta$ than in the referred references: we shall need second-order 
differentiability to establish second-order optimality conditions.

We should also mention that \cite[Lemma 2.5]{MinchenkoLeschov16} provides an arc with 
similar properties. However, the relation in (\ref{arc2}) is valid for the bigger set $J$. 
Moreover, they consider a direction in the relative interior of the critical cone, while 
here we take an arbitrary direction in the linearized cone.

A first (and immediate) consequence of Proposition \ref{prop:arc} is given next. 

\begin{corollary}
\label{corl:acq}
Suppose that RCRCQ holds at $\bar{x}\in\Omega$. Then 
$\mathcal{L}(\bar{x})\subset\mathcal{T}(\bar{x})$, which in turn implies ACQ. 
\end{corollary}

Now we derive the main result of the paper as another consequence of 
Proposition \ref{prop:arc}, with a quite elementary proof.  

\begin{theorem}
\label{th:ssonc}
Let $x^*\in\Omega$ be a local minimizer of the problem (\ref{problem}) satisfying RCRCQ. 
Then SSONC (Definition \ref{def:ssonc}) holds at $x^*$.
\end{theorem}
\beginproof
Since $\mathcal{C}^S(x^*)\subset\mathcal{L}(x^*)$, Proposition \ref{prop:arc} 
ensures the existence of a twice continuously differentiable arc $\zeta:(-\delta,\delta)\to\R^n$ 
such that $\zeta(0)=x^*$, $\zeta'(0)=d$, $\zeta(t)\in\Omega$ for all $t\in[0,\delta)$ and 
$g_j(\zeta(t))=0$ for all $t\in(-\delta,\delta)$ and $j\in I_g^+(x^*)$ (note that 
$I_g^+(x^*)\subset J$). Define $\beta:(-\delta,\delta)\to\R$ by $\beta(t)=f(\zeta(t))$. Since 
$x^*$ satisfies the KKT conditions, let $(\mu,\lambda)\in\R^m_+\times\R^p$ be an arbitrary multiplier vector associated with it. Then 
\begin{equation}
\label{kkt}
\nabla f(x^*)+\sum_{j\in I_g^+(x^*)}\mu_j\nabla g_j(x^*)+\sum_{i=1}^p\lambda_i\nabla h_i(x^*)=0
\end{equation}
and consequently $\beta'(0)=\nabla f(x^*)^Td=0$. Moreover, $x^*$ is a local minimizer of \eqref{problem} 
and hence $\beta(0)\leq\beta(t)$ for all $t\geq 0$ sufficiently small. Therefore, 
\begin{equation}
\label{beta2}
d^T\nabla^2f(x^*)d+\nabla f(x^*)^T\zeta''(0)=\beta''(0)\geq 0.
\end{equation}
On the other hand, we have 
\begin{equation}
\label{gj}
d^T\nabla^2g_j(x^*)d+\nabla g_j(x^*)^T\zeta''(0)=(g_j\circ\zeta)''(0)=0.
\end{equation}
for $j\in I_g^+(x^*)$ and
\begin{equation}
\label{h}
d^T\nabla^2h_i(x^*)d+\nabla h_i(x^*)^T\zeta''(0)=(h_i\circ\zeta)''(0)=0.
\end{equation}
for $i=1,\ldots,p$. 
Thus, multiplying (\ref{gj}) by $\mu_j$, (\ref{h}) by $\lambda_i$ and summing the resulting 
expressions over $j\in I_g^+(x^*)$ and $i=1,\ldots,p$ together with (\ref{beta2}), 
we obtain the desired result in view of (\ref{kkt}).
\endproof

It should be noted that, although RCRCQ is weaker than LICQ and CRCQ, it is not a 
necessary assumption for strong second-order condition, as we see in the next example. 

\begin{example}
\label{ex:rcrcq_not_weakest}
Consider the two dimensional problem 
$$
\begin{array}{cl}
\displaystyle\mathop{\rm minimize }  & f(x)=x_2  \\
{\rm subject\ to } & g_1(x)=x_1^2-x_2\leq 0, \\
& g_2(x)=-x_2\leq 0.
\end{array}
$$
We have that $x^*=\left(\begin{array}{c} 0 \\ 0 \end{array}\right)$ is the global 
solution of the problem. Moreover, since 
$$
\nabla g_1(x^*)=\nabla g_2(x^*)=\left(\begin{array}{c} 0 \\ -1 \end{array}\right),
$$
MFCQ holds at $x^*$. Thus, there exist Lagrange multipliers and they must satisfy 
$$
\left(\begin{array}{c} 0 \\ 1-\mu_1-\mu_2 \end{array}\right)=
\nabla f(x^*)+\mu_1\nabla g_1(x^*)+\mu_2\nabla g_2(x^*)=0,
$$
which means that $\mu_1+\mu_2=1$. The critical cone is 
$$
\mathcal{C}^S(x^*)=\{d\in\R^2\mid d_2=0 \}
$$
and the Lagrangian Hessian at $x^*$ is 
$$
\nabla^2f(x^*)+\mu_1\nabla^2g_1(x^*)+\mu_2\nabla^2g_2(x^*)=2
\left(\begin{array}{cc} \mu_1 & 0 \\ 0 & 0 \end{array}\right),
$$
so that SSONC is satisfied at this point. On the other hand, the vectors 
$\nabla g_1(x)=\left(\begin{array}{c} 2x_1 \\ -1 \end{array}\right)$ and 
$\nabla g_2(x)=\left(\begin{array}{c} 0 \\ -1 \end{array}\right)$
do not have the same rank for every $x$ in a neighborhood of $x^*$, that is, 
RCRCQ does not hold at $x^*$.
\end{example}

In view of the example above, it is natural to ask whether there is a constraint 
qualification, weaker than RCRCQ, that implies SSONC. 
In \cite[Theorem 3.2]{AndreaniBehlingHaeserSilva17} 
the authors establish SSONC under an Abadie-type assumption. Despite their hypothesis is 
weaker than RCRCQ, it is not a constraint qualification. So, the existence of Lagrange 
multipliers is not guaranteed by such an assumption.
On the other hand, the {\em second-order cone-continuity property} 
(CCP2) \cite{AndreaniHaeserRamosSilva17} is a constraint qualification implied by 
RCRCQ. However, it ensures only the {\em weak second-order optimality condition} 
(WSONC) \cite[Theorem 4.2]{AndreaniHaeserRamosSilva17}. In this reference the authors 
also define {\em strong CCP2} (SCCP2) and prove that SSONC is valid under this 
CQ \cite[Theorem 4.16]{AndreaniHaeserRamosSilva17}, but RCRCQ and SCCP2 are independent 
of each other. Both CCP2 and SCCP2 are constraint qualifications associated to a 
second-order sequential optimality condition (AKKT2) related to the convergence of 
optimization algorithms. 
Indeed, note that Example \ref{ex:rcrcq}, together with the relations among all established CQs, 
allows us to conclude that RCRCQ is the weakest constraint qualification that ensures 
strong second-order necessary condition.

\section{Conclusions}
\label{sec:concl}
In this short note we have discussed and established the first- and second-order 
optimality conditions for nonlinear programming, assuming the relaxed constant rank 
constraint qualification, through a quite elementary approach. 
The only prerequisite used here to establish auxiliary results is the classical and 
well-known inverse function theorem. We also have generalized results 
of \cite{AndreaniEchagueSchuverdt10,Janin84} assuming RCRCQ instead of CRCQ and 
presented examples and discussions that allow us to conclude that RCRCQ is the weakest CQ 
ensuring SSONC.

% \begin{acknowledgements}
% This work was partially supported by ...
% \end{acknowledgements}

% \bibliographystyle{spmpsci}
% % \bibliography{references} 
% \bibliography{/home/ademir/Ademir/bibfiles/references} 

\end{document}